\documentclass[11pt]{article}
\usepackage{}
\usepackage{mathrsfs}
\usepackage{amsfonts}

\usepackage{CJK}
\usepackage{amsfonts,amssymb,amsmath,mathrsfs}

\usepackage{color,latexsym,amsfonts}
 \newtheorem{theorem}{Theorem}[section]
 \newtheorem{lemma}[theorem]{Lemma}
 \newtheorem{corollary}[theorem]{Corollary}
 \newtheorem{proposition}[theorem]{Proposition}
 
 \newtheorem{Definition}[theorem]{Definition}
 \newtheorem{remark}[theorem]{Remark}
 \newtheorem{condition}[theorem]{Condition}

 \def\blemma{\begin{lemma}\sl{}\def\elemma{\end{lemma}}}
 \def\btheorem{\begin{theorem}\sl{}\def\etheorem{\end{theorem}}}
 \def\bcorollary{\begin{corollary}\sl{}\def\ecorollary{\end{corollary}}}
 
 \def\bproposition{\begin{proposition}\sl{}\def\eproposition{\end{proposition}}}
 \def\bremark{\begin{remark}\sl{}\def\eremark{\end{remark}}}

 \def\beqlb{\begin{eqnarray}}\def\eeqlb{\end{eqnarray}}
 \def\beqnn{\begin{eqnarray*}}\def\eeqnn{\end{eqnarray*}}

 \def\<{\langle}\def\>{\rangle}

 \def\eqref#1{{\rm(\ref{#1})}}

\def\d{\textup{d}}

\def\e{\textup{e}}

\def\fin{\hfill$\square$}

\def\newdot{{\kern.8pt\cdot\kern.8pt}}

\def\R{\mathbb{R}}
\def\E{\mathbb{E}}
\def\P{\mathbb{P}}

\def\<{\langle}
\def\>{\rangle}
\def\Proof.{\noindent{\bf Proof.}}

\begin{document}

\begin{titlepage}
\title{\bf Density estimates for the solutions of backward stochastic differential equations driven by Gaussian processes}
\author{{Xiliang Fan$^\dag$ and Jiang-Lun Wu$^\ddag$}\\
{\small $^\dag$School of Mathematics and Statistics, Anhui Normal University, Wuhu 241003, China}\\
{\small $^\ddag$Department of Mathematics, Computational Foundry, Swansea University, Swansea SA1 8EN, UK} \\
{\small \sf{fanxiliang0515@163.com} (X. Fan), \ \sf{j.l.wu@swansea.ac.uk (J.-L. Wu)}}}
\end{titlepage}
\maketitle

{\narrower{\narrower

\noindent{\bf Abstract.} The aim of this paper is twofold. Firstly, we derive upper and lower non-Gaussian bounds for the densities of the marginal laws of the solutions to backward stochastic differential equations (BSDEs) driven by fractional Brownian motions. Our arguments  consist of utilising a relationship between fractional BSDEs and quasilinear partial differential equations of mixed type, together with the profound Nourdin-Viens formula. In the linear case,
upper and lower Gaussian bounds for the densities and the tail probabilities of solutions
are obtained with simple arguments by their explicit expressions in terms of the quasi-conditional expectation. Secondly, we are concerned with Gaussian estimates for the densities of a BSDE driven by a Gaussian process in the manner that the solution can be established via an auxiliary BSDE driven by a Brownian motion. Using the transfer theorem we succeed in deriving  Gaussian estimates for the solutions.}


\bigskip



\section{Introduction}

\setcounter{equation}{0}

The problem of density estimates for solutions of stochastic equations has been extensively studied in recent years, see e.g. the monograph \cite{Nualart06a} and references therein.
Remarkably, the celebrated Bouleau-Hirsch criterion (see \cite[Theorem 2.1.2]{Nualart06a}) provides a sufficient condition for a random variable possessing a density. Moreover,
in \cite{Nourdin&Viens09a}, Nourdin and Viens derive a formula for a Malliavin differentiable random variable to admit a density with lower and upper Gaussian estimates. These results have been further extended and applied to solutions of stochastic differential equations (SDEs) and stochastic partial differential equations (SPDEs), among which let us just mention, for examples, the works by Debussche and Romito \cite{Debussche&Romito14a}, Delarue, Menozzi and Nualart \cite{Delarue&Menozzi&Nualart14a}, Millet and Sanz-Sol\'{e} \cite{Millet&Sanz-Sole99a}, Mueller and Nualart \cite{Mueller&Nualar08a}, Nualart and Quer-Sardanyons
\cite{Nualart&Quer-Sardanyons12a}, and the references therein.

On the other hand, in the seminal paper \cite{Pardoux&Peng90a} Pardoux and Peng initiated the theory of nonlinear backward stochastic differential equations (BSDEs), which is of increasing importance in stochastic control and mathematical finance (see, e.g., \cite{Karoui&Peng&Quenez97a} and most recently \cite{JianfengZhang}). This class of equations
is of the following form
\beqlb\label{1.1}
y_t=\xi+\int_t^Tf(s,y_s,z_s)\d s-\int_t^Tz_s\d B_s, \ \ t\in[0,T]
\eeqlb
on a given filtered probability space $(\Omega,\mathcal{F},\mathbb{P};\{\mathcal{F}_t\}_{t\in[0,T]})$, where $T>0$ is arbitrarily fixed, $\xi$ is a $\mathcal{F}_T$-measurable random variable, the generator $f:[0,T]\times\mathbb{R}\times\mathbb{R}\to\mathbb{R}$ is a jointly measurable map, and $B=(B_t)_{t\ge0}$ is a Brownian motion adapted to $\{\mathcal{F}_t\}_{t\in[0,T]}$£¬ or simply $\{\mathcal{F}_t\}_{t\in[0,T]}$ is taken as the natural filtration of $B$. Recall that
a solution to the BSDE \eqref{1.1} is a pair of predictable processes $(y,z)$ with suitable integrability conditions such that \eqref{1.1} holds $\P$-a.s.. To date, there is a wealth of existence and uniqueness results under various assumptions on the generators $f$ including, for instnace, the cases of Lipschitz or (super-)quadratic growth \cite{Delbaen&Hu&Bao10a,Karoui&Peng&Quenez97a,Kobylanski00a,Pardoux&Peng90a,Richou12a}.
When dealing with applications such as the numerical approximation of the solutions,
one needs to investigate the existence and regularity of densities for the marginal laws of $(y,z)$. As far as we know, there are comparably only a few works to study this problem. The first results have been derived by Antonelli and Kohatsu-Higa \cite{Antonelli&Kohatsu-Higa05a}, in which they study the existence and the estimates of the density for $y_t$ at a fixed time $t\in[0,T]$ via the Bouleau-Hirsch criterion. Then, based on the Nourdin-Viens formula, Aboura and Bourguin \cite{Aboura&Bourguin13a} have proved the existence of the density for $z_t$ under the condition that the generator $f$ is linear with respect to its $z$ variable, and further obtained the estimates on the densities of the laws of $y_t$ and $z_t$. Recently, Mastrolia, Possama\"{i} and R\'{e}veillac in \cite{Mastrolia&Possama&Reveillac16a} have studied the existence of densities for marginal laws of the solution $(y,z)$ to \eqref{1.1} with a quadratic growth generator, and derived the estimates on these densities. Afterwards, Mastrolia \cite{Mastrolia16a} has extended the results to the case of non-Markovian BSDEs.

One of the main objective of the present paper concerns the problem of density estimates for the following BSDE
\beqlb\label{1.2}
y_t=h(\eta_T)+\int_t^Tf(s,\eta_s,y_s,z_s)\d s-\int_t^Tz_s\d B^H_s, \ \ t\in[0,T],
\eeqlb
where $\eta_t=\eta_0+\int_0^tb_s\d s+\int_0^t\sigma_s\d B^H_s$ with $\eta_0, b_s$ and $\sigma_s$ being respectively a constant and deterministic functions, and $B^H$ is a fractional Brownian motion with Hurst parameter $H\in(0,1)$,
the stochastic integral is the divergence-type integral.
Precise assumptions on the (deterministic and joint measurable) generator $f:[0,T]\times\R\times\R\times\R\rightarrow\R$ and $h:\R\rightarrow\R$
will be specified in later sections. Let us recall that, $B^H$ with Hurst parameter $H\in(0,1)$ is a center Gaussian process with covariance
\beqnn
 R_H(t,s):=\E\left(B^H_tB^H_s\right)=\frac{1}{2}\left(t^{2H}+s^{2H}-|t-s|^{2H}\right),\, t,s\in[0,T].
\eeqnn
This implies that for each $p\geq 1$, there holds $\E(|B_t^H-B_s^H|^p)=C(p)|t-s|^{pH}$.
Then $B^H$ is $(H-\epsilon)$-order H\"{o}lder continuous for any $\epsilon>0$ and is an $H$-self similar process. This, together with the fact that $B^{1/2}$ is a Brownian motion, converts fractional Brownian motion into a natural generalization of Brownian motion
and leads to many applications in modelling physical phenomena and finance behaviours.

We mention that there are several papers concerning the existence and uniqueness results of solutions for \eqref{1.2}.
Biagini, Hu, \O ksendal and Sulem in \cite{Biagini&Hu&Oksendal&Sulem02a} first studied linear fractional BSDE with $H\in(1/2,1)$ which are based on fractional Clark-Ocone formula and the Girsanov transformation. In the spirit of the four step scheme introduced by Ma, Protter and
Yong \cite{Ma&Protter&Yong94a} for BSDEs perturbed by a standard Brownian motion,
Bender \cite{Bender05a} constructed an explicit solutions for a kind of linear fractional BSDEs with $H\in(0,1)$ via the solution of some PDE and fractional It\^{o} formula.
In the case of nonlinear fractional BSDEs with $H\in(1/2,1)$,
Hu and Peng \cite{Hu&Peng09a} first proved the existence and uniqueness of the solution through the notion of quasi-conditional expectation introduced in \cite{Hu&Oksendal03a}.
Then, based on \cite{Hu&Peng09a}, Maticiuc and Nie in \cite{Maticiuc&Nie15a} made some improvements of analysis and extended to the case of fractional backward stochastic variational inequalities, and Fei, Xia and Zhang in \cite{Fei&Xia&Zhang13a} generalized the investigation to BSDEs driven by both standard and fractional Brownian motions.
In a multivariate setting where each of the components is an independent fractional Brownian motion $B_i^{H_i}$ with $H_i\in(1/2,1)$, Hu, Ocone and Song in \cite{Hu&cone&Song12a} solved fractional BSDEs by using their relation to PDEs, and they further derived a comparison theorem.

In the present paper, with the help of the connection between the solution to BSDE \eqref{1.2} and the solution to its associated PDE of mixed type, we shall give some sufficient conditions to ensure the existence of densities for marginal laws of the solution $(y,z)$ to BSDE \eqref{1.2}. Moreover, we will derive non-Gaussian tail estimates of densities.
To the best of our knowledge these kind of estimates for BSDE \eqref{1.2} are not available in the existing literature.
When
\beqnn
f(s,\eta_s,y_s,z_s)=\alpha_s+\beta_s y_s+\gamma_t z_s,
\eeqnn
i.e., BSDE \eqref{1.2} is linear,
the Gaussian bounds for the densities and the tail probabilities of solutions
will be derived with a direct and simpler arguments by their explicit expressions in terms of the quasi-conditional expectation.

Our paper is also dedicated to obtaining Gaussian bounds for the densities of the solution to BSDE
\beqlb\label{1.3}
y_t=h(X_T)+\int_t^Tf(s,X_s,y_s,z_s)\d V(t)-\int_t^Tz_s\d^\diamond X_s,
\eeqlb
where $X$ is a centered Gaussian process with a strictly increasing continuous variance function $V(t)=\mathrm{Var}X_t$, the stochastic integral is the Wick-It\^{o} integral defined via the $S$-transformation and the Wick product. When $X$ is a Brownian motion, the above Wick-It\^{o} integral coincides with the classical It\^{o} integral. In \cite{Bender} Bender shows the existence and uniqueness results and then obtains a strict comparison theorem for BSDE \eqref{1.3} using the transfer theorem which can transfer the concerned problems to an auxiliary BSDE driven by a Brownian motion. In \cite{Bender} the author also compares this type of equations with other BSDEs driven by Gaussian non-semimartingales, especially BSDEs driven by fractional Brownian motion $B^H$ with $H\in(1/2,1)$. The final objective of the present paper is to deepen the investigation of BSDE \eqref{1.3}. We study the Gaussian bounds for marginal laws of the solution $(y,z)$ to BSDE \eqref{1.3} via the transfer theorem.

The rest of our paper is structured as follows. Section 2, the next section, presents some basic elements of stochastic calculus with respect to fractional Brownian motion which are needed in later sections. We investigate the non-Gaussian bounds for the densities of the nonlinear fractional BSDEs in Section 3.
Section 4 is devoted to the derivation of the Gaussian bounds for the densities and the tail probabilities of linear fractional BSDEs. Finally in Section 5, we provide the Gaussian bounds for the densities of BSDEs driven by Gaussian processes.

\section{Preliminaries}

\setcounter{equation}{0}

In this section, we shall give some basic elements of stochastic calculus with respect to fractional Brownian motion. For a deeper and detailed discussion, we refer the reader to \cite{Alos&Mazet&Nualart01a,Biagini&Hu08a,Decreusefond&Ustunel98a} and \cite{Nualart06a}.

Let $\Omega$ be the canonical probability space $C_0([0,T],\R)$, i.e., the Banach space of continuous functions on $[0,T]$ vanishing at time $0$, equipped with the supremum norm, and $\mathcal{F}$ is taken to be the Borel $\sigma$-algebra. Let $\P$ be the unique probability measure on $\Omega$ such that the canonical process $B^H=(B^H_t)_{t\in[0,T]}$ is a fractional Brownian motion with Hurst parameter $H\in(0,1)$.

Let $\mathscr{E}$ be the space of step functions on $[0,T]$, and $\mathcal {H}$  the closure of $\mathscr{E}$ with respect to the following scalar product determined by the covariance $R_H$ of $B^H$
\beqnn
\langle I_{[0,t]},I_{[0,s]}\rangle_\mathcal {H}=R_H(t,s).
\eeqnn
By the bounded linear transformation theorem, 
the mapping $I_{[0,t]}\mapsto B_t^H$ can be extended to an isometry between $\mathcal {H}$ and the Gaussian space $\mathcal {H}_1$ associated with $B^H$.
Denote this isometry by $\phi\mapsto B^H(\phi)$.
When $H\in(1/2,1)$ it can be shown that $L^{1/H}([0,T])\subset\mathcal{H}$,
and when $H\in(0,1/2)$ there holds $\mathcal{H}\subset L^2([0,T])$.
For $H\in(1/2,1)$, we shall use the following representation of the inner product in $\mathcal{H}$:
\beqlb\label{2.1.1}
\langle\phi,\psi\rangle_\mathcal{H}=C_H\int_0^T\int_0^T\phi_u\psi_v|u-v|^{2H-2}\d u\d v,
\eeqlb
where $C_H:=H(2H-1)$.

Let $\mathcal {S}$ denote the totality of smooth and cylindrical random variables of the form
\beqnn
F=f(B^H(\phi_1),\cdot\cdot\cdot,B^H(\phi_n)),
\eeqnn
where $n\geq 1, f\in C_b^\infty(\mathbb{R}^n)$, the set of $f$ and all its partial derivatives are bounded, $\phi_i\in\mathcal{H}, 1\leq i\leq n$. The Malliavin derivative of $F$, denoted by $D^HF$, is defined as the $\mathcal {H}$-valued random variable
\beqnn
D^HF=\sum_{i=1}^n\frac{\partial f}{\partial x_i}(B^H(\phi_1),\cdot\cdot\cdot,B^H(\phi_n))\phi_i.
\eeqnn
For any $p\geq 1$, we define the Sobolev space $\mathbb{D}_H^{1,p}$ as the completion of $\mathcal {S}$ with respect to the norm
\beqnn
\|F\|_{1,p}^p=\mathbb{E}|F|^p+\mathbb{E}\|D^HF\|^p_{\mathcal {H}}.
\eeqnn

Next, let $F$ be in $\mathbb{D}_H^{1,2}$ and write $D^HF:=\Phi_F(B^H)$,
where $\Phi_F:\R^\mathcal{H}\rightarrow\mathcal{H}$ is a measurable mapping.
Set
\beqlb\label{2.2.1}
g_F(x)=\int_0^\infty\e^{-\theta}\E\left[\E'\left[\langle\Phi_F(B^H),
\widetilde{\Phi_F^\theta}(B^H)\rangle_\mathcal{H}\right]|F-\E F=x\right]\d \theta,\ x\in\R
\eeqlb
where $\widetilde{\Phi_F^\theta}(B^H):=\Phi_F(\e^{-\theta}B^H+\sqrt{1-\e^{-2\theta}}B'^H)$
with $B'^H$ an independent copy of $B^H$ such that $B^H$ and $B'^H$ are defined on the product probability space
$(\Omega\times\Omega',\mathcal{F}\times\mathcal{F'},\P\times\P')$.
We recall the following result, cf. \cite[Theorem 3.1 and Proposition 3.7]{Nourdin&Viens09a},
which presents a criterion for a Malliavin differentiable random variable to have a density with Gaussian bounds
based on the above function $g$.

\bproposition\label{Theor2.2}
The law of $F$ has a density $\rho_F$ with respect to the Lebesgue measure if and only if $g_F(F-\E F)>0$ a.s.
In this case, $\mathrm{Supp}(\rho_F)$ is a closed interval of $\R$ and for all $z\in\mathrm{Supp}(\rho_F)$, there holds
\beqnn
\rho_F(z)=\frac{\E |F-\E F|}{2g_F(z-\E F)}\exp\left(-\int_0^{z-\E F}\frac{u\d u}{g_F(u)}\right).
\eeqnn
Furthermore, if there exist constants $c_1,c_2>0$ such that
\beqnn
c_1\leq g_F(x)\leq c_2,\ \ \ \P-\mathrm{a.s.},
\eeqnn
then the law of $F$ has a density $\rho$ satisfying, for almost all $x\in\R$,
\beqnn
\frac{\E|F-\E F|}{2c_2}\exp\left(-\frac{(x-\E F)^2}{2c_1}\right)\leq\rho(x)\leq\frac{\E|F-\E F|}{2c_1}\exp\left(-\frac{(x-\E F)^2}{2c_2}\right).
\eeqnn
\eproposition

Besides, by \cite[Theorem 4.1]{Nourdin&Viens09a} (see also \cite[Proposition 2.2]{Nguyen&Privault&Torrisi}) we then have the following tail estimates.
\bproposition\label{Proposition 2.1}
Let $F\in\mathbb{D}_H^{1,2}$ with $\E F=0$.
If $0<g_F(x)\leq a_1x+a_2,$ a.s. for some $a_1\geq0$ and $a_2>0$,
then
\beqnn
\P(F\geq x)\leq\exp\left(-\frac{x^2}{2a_1 x+2a_2}\right)\ \
\mathrm{and} \ \ \ \P(F\leq-x)\leq\exp\left(-\frac{x^2}{2a_2}\right),\ \ x>0.
\eeqnn
\eproposition

\section{BSDEs driven by fractional Brownian motions}

\setcounter{equation}{0}

The objective of this section is to study the non-Gaussian densities estimates for the solution of the following BSDE driven by fractional Brownian motion
\begin{equation}\label{3.1}
\left\{
\begin{array}{ll}
\d y_t=f(t,\eta_t,y_t,z_t)\d t-z_t\d B^H_t,\\
y_T=h(\eta_T),
\end{array} \right.
\end{equation}
with
\beqnn
\eta_t:=\eta_0+\int_0^tb_s\d s+\int_0^t\sigma_s\d B^H_s.
\eeqnn
Here $\eta_0$ is a given constant, $b$ and $\sigma$ are bounded deterministic functions such that $\sigma(t)\neq 0$ for all $t\in[0,T], H\in(\frac{1}{2},1)$.
A pair of $\mathcal{F}_t$-adapted stochastic
processes $(y,z)$ is called a solution to the equation \eqref{3.1} if
\beqnn
y_t=h(\eta_T)+\int_t^Tf(s,\eta_s,y_s,z_s)\d s-\int_t^Tz_s\d B^H_s, \ \ t\in[0,T].
\eeqnn
We begin with the assumption (H1) below
\begin{itemize}
\item[(i)] $f:[0,T]\times\R^3\rightarrow\R$ is differentiable with respect to the third component and there exists a nonnegative constant $K$ such that,
for all $t\in[0,T], x, y_1,y_2,z_1,z_2\in\R$,
\beqnn
|f(t,x,y_1,z_1)-f(t,x,y_2,z_2)|+|f_y(t,x,y_1,z_1)-f_y(t,x,y_2,z_2)|\leq K(|y_1-y_2|+|z_1-z_2|).
\eeqnn

\par

\item[(ii)]$h:\R\rightarrow\R$ is continuously differentiable and of polynomial growth.
\end{itemize}

Due to \cite[Theorem 3.4]{Hu&cone&Song12a}, the condition (H1) ensures that there exists a unique solution $(y,z)$ to BSDE \eqref{3.1},
which is given by $y_t=u(t,\eta_t)$ and $v_t=\sigma_tu_x(t,\eta_t)$ via some deterministic function $u:[0,T]\times\R\rightarrow\R$, where $u_x(t,x):=\frac{\partial}{\partial x}u(t,x)$.
The argument in \cite{Hu&cone&Song12a} is based on a connection between this equation and a quasilinear PDE of mixed type.
In the remaining part of this section, we assume (H1) holds and moreover the unique solution is of the above form.
We aim to show the non-Gaussian densities estimates for the marginal laws of $(y,z)$ at a fixed time $t\in(0,T)$.
To this end, we let
\beqnn
\varrho_t=C_H\int_0^t\int_0^t\sigma_u\sigma_v|u-v|^{2H-2}\d u\d v,\ \ \ p_{\varrho_t}(x)=\frac{1}{\sqrt{2\pi\varrho_t}}\e^{-\frac{x^2}{2\varrho_t}}
\eeqnn
and for each $h\in L^0(\R)$, define
\beqnn
\overline{h}:=\inf\{\gamma>0:\limsup\limits_{|x|\rightarrow+\infty}\frac{|h(x)|}{|x|^\gamma}<\infty\},\ \ \
\underline{h}:=\inf\{\gamma>0:\liminf\limits_{|x|\rightarrow+\infty}\frac{|h(x)|}{|x|^\gamma}<\infty\}.
\eeqnn
Clearly, the above $\overline{h}$ and $\underline{h}$ can describe the asymptotic growth  of $h$ in the neighborhood of $+\infty$ and $-\infty$.

Our main result of this section reads as follows
\btheorem\label{Theorem4.1}
Let $t\in(0,T]$. Suppose that $0<\underline{u(t,\cdot)}<+\infty,\overline{u_x(t,\cdot)}<+\infty$ and there exist positive constants $L,\lambda$ satisfying
$u_x(t,\cdot)\geq\frac{1}{L(1+|\cdot|^\lambda)}$.
Then, the law of $y_t=u(t,\eta_t)$ has a density $\rho_{y_t}$,
and for any $\epsilon,\delta>0$ there exist positive constants $C_{\epsilon,t}$ and $C_{\delta,t}$ depending on $\epsilon,t$ and $\delta,t$, respectively, such that
\beqlb\label{add-Theorem4.1-1}
&&\frac{\E |y_t-\E y_t|}{2C(\epsilon,\delta)\varrho_t\left(1+|z|^{\bar{\epsilon}\bar{\delta}}\right)
\left(1+|z|^{\bar{\epsilon}\bar{\delta}}+\varrho_t^{\frac{\bar{\epsilon}}{2}}\right)}
\exp\left(-\frac{1}{\tilde{C}(\delta)\varrho_t}\int_0^{z-\E y_t}u\left(1+|u+\E{y_t}|^{2\lambda\bar{\delta}}\right)\d u\right)\leq
\rho_{y_t}(z)\nonumber\\
&\leq&\frac{\E |y_t-\E y_t|}{2\tilde{C}(\delta)\varrho_t}
\left(1+|z|^{2\lambda\bar{\delta}}\right)\exp\left(-\frac{1}{C(\epsilon,\delta)\varrho_t}\int_0^{z-\E y_t}
\frac{u\d u}{\left(1+|u+\E{y_t}|^{\bar{\epsilon}\bar{\delta}}\right)\left(1+|u+\E{y_t}|^{\bar{\epsilon}\bar{\delta}}+\varrho_t^{\frac{\bar{\epsilon}}{2}}\right)}\right),\nonumber\\
\eeqlb
where
\beqnn
C(\epsilon,\delta)&:=&\sup\limits_{t\in[0,T]}\Bigg[C_{\epsilon,t}^2\left(1+C_{\delta,t}^{\bar{\epsilon}}2^{(\bar{\epsilon}-1)_+}\right)\nonumber\\
&&\ \ \ \ \ \ \ \ \ \cdot\left(\left(1+3^{(\bar{\epsilon}-1)_+}\left|\eta_0+\int_0^tb_s\d s\right|^{\bar{\epsilon}}
+\frac{C_{\delta,t}^{\bar{\epsilon}}6^{(\bar{\epsilon}-1)_+}}{1+\bar{\epsilon}}\right)
\vee\left(\frac{\Gamma\left(\frac{1+\bar{\epsilon}}{2}\right)2^{\frac{\bar{\epsilon}}{2}}}{\sqrt{\pi}}\right)\right)\Bigg]
\eeqnn
and
\beqnn
\tilde{C}(\delta):=\frac{1}{3^{(\lambda-1)_+}2L^2(1+2^{(2\lambda-1)_+}C_{\delta,t}^{2\lambda})}
\cdot\sup\limits_{t\in[0,T]}\int_\R\frac{p_{\varrho_t}(z)}{1+\left|\eta_0+\int_0^tb_s\d s\right|^\lambda+|z|^\lambda}\d z
\eeqnn
with $\bar{\epsilon}:=\overline{u_x(t,\cdot)}+\epsilon$ and $\bar{\delta}:=\overline{u^{-1}(t,\cdot)}+\delta$.
\etheorem

\emph{Proof.} By the fact that $u(t,\cdot)$ is continuous and increasing,
we easily verify that $y_t$ has a density $\rho_{y_t}$.
In order to prove \eqref{add-Theorem4.1-1},
we rely heavily on Proposition \ref{Theor2.2}.

Notice first that for $0<u\leq t\leq T$, we have $D^H_uy_t=u_x(t,\eta_t)\sigma_u$.
Then, it follows that
\beqnn
\Phi_{y_t}(B^H)=u_x(t,\eta_t)\sigma_\cdot
\eeqnn
and
\beqnn
\widetilde{\Phi_{y_t}^\theta}(B^H)&=&\Phi_{y_t}\left(e^{-\theta}B^H+\sqrt{1-\e^{-2\theta}}B'^H\right)\\
&=&u_x\left(t,\left(1-e^{-\theta}\right)\left(\eta_0+\int_0^tb_s\d s\right)+e^{-\theta}\eta_t+\sqrt{1-\e^{-2\theta}}\int_0^t\sigma_s\d B'^H_s\right)\sigma_\cdot.
\eeqnn
Hence, using \eqref{2.2.1} and \eqref{2.1.1}, we deduce that for $y\in\R$,
\beqlb\label{Theorem4.1-1}
&&g_{y_t}(y)\nonumber\\
&=&\int_0^\infty\e^{-\theta}\E\left[\E'\left[\langle\Phi_{y_t}(B^H),\widetilde{\Phi_{y_t}^\theta}(B^H)\rangle_\mathcal{H}\right]\Big|{y_t}-\E {y_t}=y\right]\d \theta\nonumber\\
&=&\varrho_t
\int_0^\infty\e^{-\theta}\E\Big[u_x(t,\eta_t)\nonumber\\
&&\E'u_x\left(t,\left(1-e^{-\theta}\right)\left(\eta_0+\int_0^tb_s\d s\right)+e^{-\theta}\eta_t+\sqrt{1-\e^{-2\theta}}\int_0^t\sigma_s\d B'^H_s\right)
\Big|\eta_t=u^{-1}(t,y+\E{y_t})\Big]\d \theta\nonumber\\
&=&\varrho_t u_x(t,u^{-1}(t,y+\E{y_t}))\nonumber\\
&&\cdot\int_0^\infty\e^{-\theta}
\int_\R u_x\left(t,\left(1-e^{-\theta}\right)\left(\eta_0+\int_0^tb_s\d s\right)+e^{-\theta}u^{-1}(t,y+\E{y_t})+\sqrt{1-\e^{-2\theta}}z\right)
p_{\varrho_t}(z)\d z\d\theta.\nonumber\\
\eeqlb
With the relation \eqref{Theorem4.1-1} in hand, we shall follow the strategy designed in \cite{Mastrolia&Possama&Reveillac16a} to obtain upper and lower bounds for $g_{y_t}$.

\textsl{Upper bound}.
Since $\underline{u(t,\cdot)}\in(0,\infty)$, we know from \cite[Lemma 5.2]{Mastrolia&Possama&Reveillac16a} that $\overline{u^{-1}(t,\cdot)}\in(0,\infty)$.
Taking into account the definitions of $\overline{u_x(t,\cdot)}$ and $\overline{u^{-1}(t,\cdot)}$, we have for each $\epsilon>0$
\beqlb\label{addTheorem4.1-1}
(0<)u_x(t,z)\leq C_{\epsilon,t}(1+|z|^{\overline{u_x(t,\cdot)}+\epsilon}),\ \ \forall z\in\R,
\eeqlb
and for each $\delta>0$
\beqlb\label{addTheorem4.1-2}
|u^{-1}(t,z)|\leq C_{\delta,t}(1+|z|^{\overline{u^{-1}(t,\cdot)}+\delta}),\ \ \forall z\in\R,
\eeqlb
where $C_{\epsilon,t}$ and $C_{\delta,t}$ are both constants depending on $\epsilon,t$ and $\delta,t$ respectively.\\
For the convenience of the notation, we let $\bar{\epsilon}=\overline{u_x(t,\cdot)}+\epsilon$ and $\bar{\delta}=\overline{u^{-1}(t,\cdot)}+\delta$.
Then plugging the inequalities \eqref{addTheorem4.1-1} and \eqref{addTheorem4.1-2} and resorting to the $C_r$-inequality, we have
\beqlb\label{Theorem4.1-2}
&&g_{y_t}(y)\nonumber\\
&\leq&\varrho_t C_{\epsilon,t}^2\left(1+|u^{-1}(t,y+\E{y_t})|^{\bar{\epsilon}}\right)\nonumber\\
&&\cdot\int_0^\infty\e^{-\theta}
\int_\R\left(1+\left|\left(1-e^{-\theta}\right)\left(\eta_0+\int_0^tb_s\d s\right)+e^{-\theta}u^{-1}(t,y+\E{y_t})+\sqrt{1-\e^{-2\theta}}z\right|^{\bar{\epsilon}}\right)p_{\varrho_t}(z)\d z\d\theta
\nonumber\\
&\leq&\varrho_t C_{\epsilon,t}^2\left(1+|u^{-1}(t,y+\E{y_t})|^{\bar{\epsilon}}\right)\nonumber\\
&&\cdot\int_0^\infty\e^{-\theta}
\int_\R\left(1+3^{(\bar{\epsilon}-1)_+}\left[\left|\eta_0+\int_0^tb_s\d s\right|^{\bar{\epsilon}}+e^{-\bar{\epsilon}\theta}
|u^{-1}(t,y+\E{y_t})|^{\bar{\epsilon}}+|z|^{\bar{\epsilon}}\right]\right)p_{\varrho_t}(z)\d z\d\theta
\nonumber\\
&=&\varrho_t C_{\epsilon,t}^2\left(1+|u^{-1}(t,y+\E{y_t})|^{\bar{\epsilon}}\right)\nonumber\\
&&\cdot\left(1+3^{(\bar{\epsilon}-1)_+}\left|\eta_0+\int_0^tb_s\d s\right|^{\bar{\epsilon}}
+\frac{3^{(\bar{\epsilon}-1)_+}}{1+\bar{\epsilon}}|u^{-1}(t,y+\E{y_t})|^{\bar{\epsilon}}
+\frac{\Gamma\left(\frac{1+\bar{\epsilon}}{2}\right)2^{\frac{\bar{\epsilon}}{2}}}{\sqrt{\pi}}\varrho_t^{\frac{\bar{\epsilon}}{2}}\right)\nonumber\\
&\leq&\varrho_t C_{\epsilon,t}^2
\left(1+C_{\delta,t}^{\bar{\epsilon}}2^{(\bar{\epsilon}-1)_+}\left(1+|y+\E{y_t}|^{\bar{\epsilon}\bar{\delta}}\right)\right)\nonumber\\
&&\cdot\left(1+3^{(\bar{\epsilon}-1)_+}\left|\eta_0+\int_0^tb_s\d s\right|^{\bar{\epsilon}}
+\frac{C_{\delta,t}^{\bar{\epsilon}}6^{(\bar{\epsilon}-1)_+}}{1+\bar{\epsilon}}\left(1+|y+\E{y_t}|^{\bar{\epsilon}\bar{\delta}}\right)
+\frac{\Gamma\left(\frac{1+\bar{\epsilon}}{2}\right)2^{\frac{\bar{\epsilon}}{2}}}{\sqrt{\pi}}\varrho_t^{\frac{\bar{\epsilon}}{2}}\right)\nonumber\\
&\leq&C(\epsilon,\delta)\varrho_t\left(1+|y+\E{y_t}|^{\bar{\epsilon}\bar{\delta}}\right)
\left(1+|y+\E{y_t}|^{\bar{\epsilon}\bar{\delta}}+\varrho_t^{\frac{\bar{\epsilon}}{2}}\right)
\eeqlb
where
\beqnn
C(\epsilon,\delta)&=&\sup\limits_{t\in[0,T]}\Bigg[C_{\epsilon,t}^2\left(1+C_{\delta,t}^{\bar{\epsilon}}2^{(\bar{\epsilon}-1)_+}\right)\nonumber\\
&&\ \ \ \ \ \ \ \ \ \cdot\left(\left(1+3^{(\bar{\epsilon}-1)_+}\left|\eta_0+\int_0^tb_s\d s\right|^{\bar{\epsilon}}
+\frac{C_{\delta,t}^{\bar{\epsilon}}6^{(\bar{\epsilon}-1)_+}}{1+\bar{\epsilon}}\right)
\vee\left(\frac{\Gamma\left(\frac{1+\bar{\epsilon}}{2}\right)2^{\frac{\bar{\epsilon}}{2}}}{\sqrt{\pi}}\right)\right)\Bigg].
\eeqnn

\textsl{Lower bound}.
Due to the condition on $u_x(t,\cdot)$, the $C_r$-inequality and \eqref{addTheorem4.1-2},  we obtain
\beqlb\label{Theorem4.1-3}
&&g_{y_t}(y)\nonumber\\
&\geq&\frac{\varrho_t }{L^2(1+|u^{-1}(t,y+\E{y_t})|^\lambda)}\nonumber\\
&&\cdot\int_0^\infty\e^{-\theta}\int_\R
\frac{p_{\varrho_t}(z)}{1+\left|\left(1-e^{-\theta}\right)\left(\eta_0+\int_0^tb_s\d s\right)+e^{-\theta}u^{-1}(t,y+\E{y_t})+\sqrt{1-\e^{-2\theta}}z\right|^\lambda}\d z\d\theta\nonumber\\
&\geq&\frac{\varrho_t }{L^2(1+|u^{-1}(t,y+\E{y_t})|^\lambda)}
\cdot\int_\R\frac{p_{\varrho_t}(z)}
{1+3^{(\lambda-1)_+}\left(\left|\eta_0+\int_0^tb_s\d s\right|^\lambda+\left|u^{-1}(t,y+\E{y_t})\right|^\lambda+|z|^\lambda\right)}
\d z\nonumber\\
&\geq&\frac{1}{3^{(\lambda-1)_+}L^2}
\cdot\int_\R\frac{p_{\varrho_t}(z)}
{1+\left|\eta_0+\int_0^tb_s\d s\right|^\lambda+|z|^\lambda}\d z
\cdot\frac{\varrho_t }{(1+|u^{-1}(t,y+\E{y_t})|^\lambda)^2}\nonumber\\
&\geq&\frac{1}{3^{(\lambda-1)_+}2L^2}
\cdot\int_\R\frac{p_{\varrho_t}(z)}{1+\left|\eta_0+\int_0^tb_s\d s\right|^\lambda+|z|^\lambda}\d z
\cdot\frac{\varrho_t}{1+C_{\delta,t}^{2\lambda}\left(1+|y+\E{y_t}|^{\bar{\delta}}\right)^{2\lambda}}\nonumber\\
&\geq&\tilde{C}(\delta)\frac{\varrho_t}{1+|y+\E{y_t}|^{2\lambda\bar{\delta}}},
\eeqlb
where $\tilde{C}(\delta)=\frac{1}{3^{(\lambda-1)_+}2L^2(1+2^{(2\lambda-1)_+}C_{\delta,t}^{2\lambda})}
\cdot\sup\limits_{t\in[0,T]}\int_\R\frac{p_{\varrho_t}(z)}{1+\left|\eta_0+\int_0^tb_s\d s\right|^\lambda+|z|^\lambda}\d z$.\\
Therefore, applying Proposition \ref{Theor2.2}, together with \eqref{Theorem4.1-2} and \eqref{Theorem4.1-3}, we complete the proof.
\fin

\bremark\label{Remark4.1}
Let us have a close look at the lower bound of $\rho_{y_t}$. In fact, a direct computation shows that
\beqlb\label{Corollary4.1-2}
\int_0^{z-\E y_t}u\left(1+|u+\E{y_t}|^{2\lambda\bar{\delta}}\right)\d u
&=&\frac{1}{2}|z-\E y_t|^2+\frac{1}{2(1+\lambda\bar{\delta})}\left(|z|^{2(1+\lambda\bar{\delta})}-|\E y_t|^{2(1+\lambda\bar{\delta})}\right)\nonumber\\
&&-\frac{|\E y_t|}{1+2\lambda\bar{\delta}}\left(\mathrm{sgn}(z\E y_t)|z|^{1+2\lambda\bar{\delta}}-|\E y_t|^{1+2\lambda\bar{\delta}}\right)
=:\chi(z,\E y_t).\nonumber\\
\eeqlb
Then, it yields the following
\beqnn
\rho_{y_t}(z)\geq
\frac{\E |y_t-\E y_t|}{2C(\epsilon,\delta)\varrho_t\left(1+|z|^{\bar{\epsilon}\bar{\delta}}\right)
\left(1+|z|^{\bar{\epsilon}\bar{\delta}}+\varrho_t^{\frac{\bar{\epsilon}}{2}}\right)}
\exp\left(-\frac{\chi(z,\E y_t)}{\tilde{C}(\delta)\varrho_t}\right).
\eeqnn
\eremark

As for the upper bound in Theorem \ref{Theorem4.1}, we get the following result.
\bcorollary\label{Corollary4.1}
Under the assumptions in Theorem \ref{Theorem4.1}, there exists $z_0>0$ such that
\beqnn
\rho_{y_t}(z)
\leq\frac{\E |y_t-\E y_t|}{2\tilde{C}(\delta)\varrho_t}
\left(1+|z|^{2\lambda\bar{\delta}}\right)\exp\left(-\frac{\left|z-\E y_t\right|^{2(1-\bar{\epsilon}\bar{\delta})}-\left|\mathrm{sgn}z\cdot z_0-\E y_t\right|^{2(1-\bar{\epsilon}\bar{\delta})}}{4(1-\bar{\epsilon}\bar{\delta})C(\epsilon,\delta)\varrho_t}
\right)
\eeqnn
holds for all $|z|>z_0$.
\ecorollary

\emph{Proof.}
We first observe that by Theorem \ref{Theorem4.1} our problem can be reduced to show that the following inequality
\beqlb\label{Corollary4.1-1}
&&\int_0^{z-\E y_t}\frac{u\d u}{\left(1+|u+\E{y_t}|^{\bar{\epsilon}\bar{\delta}}\right)\left(1+|u+\E{y_t}|^{\bar{\epsilon}\bar{\delta}}+\varrho_t^{\frac{\bar{\epsilon}}{2}}\right)}\nonumber\\
&\geq&\frac{1}{4(1-\bar{\epsilon}\bar{\delta})}\left(\left|z-\E y_t\right|^{2(1-\bar{\epsilon}\bar{\delta})}-\left|\mathrm{sgn}z\cdot z_0-\E y_t\right|^{2(1-\bar{\epsilon}\bar{\delta})}\right)
\eeqlb
holds for each $|z|>z_0$.
In order to prove \eqref{Corollary4.1-1}, we start by noticing that
\beqnn
\lim\limits_{u\rightarrow+\infty}\frac{u}{\left(1+|u+\E{y_t}|^{\bar{\epsilon}\bar{\delta}}\right)\left(1+|u+\E{y_t}|^{\bar{\epsilon}\bar{\delta}}+\varrho_t^{\frac{\bar{\epsilon}}{2}}\right)}
\frac{1}{\frac{u}{|u|^{2\bar{\epsilon}\bar{\delta}}}}=1,
\eeqnn
and then for some sufficiently large $u_0>0$, we have for $u\geq u_0$
\beqnn
\frac{u}{\left(1+|u+\E{y_t}|^{\bar{\epsilon}\bar{\delta}}\right)\left(1+|u+\E{y_t}|^{\bar{\epsilon}\bar{\delta}}+\varrho_t^{\frac{\bar{\epsilon}}{2}}\right)}
\geq\frac{u}{2|u|^{2\bar{\epsilon}\bar{\delta}}}.
\eeqnn
Consequently, there exists $z_0>0$ large enough such that \eqref{Corollary4.1-1} holds for any $|z|\geq z_0$.
Indeed, when $z\leq-z_0$ (the choosing of $z_0$ depends on the above argument and moreover satisfies $z_0-|\E y_t|>0$), we get
\beqnn
&&\int_0^{z-\E y_t}\frac{u\d u}{\left(1+|u+\E{y_t}|^{\bar{\epsilon}\bar{\delta}}\right)\left(1+|u+\E{y_t}|^{\bar{\epsilon}\bar{\delta}}+\varrho_t^{\frac{\bar{\epsilon}}{2}}\right)}\\
&=&\int_{z-\E y_t}^0\frac{-u\d u}{\left(1+|u+\E{y_t}|^{\bar{\epsilon}\bar{\delta}}\right)\left(1+|u+\E{y_t}|^{\bar{\epsilon}\bar{\delta}}+\varrho_t^{\frac{\bar{\epsilon}}{2}}\right)}\\
&\geq&
\int_{z-\E y_t}^{-z_0-\E y_t}\frac{-u\d u}{\left(1+|u+\E{y_t}|^{\bar{\epsilon}\bar{\delta}}\right)\left(1+|u+\E{y_t}|^{\bar{\epsilon}\bar{\delta}}+\varrho_t^{\frac{\bar{\epsilon}}{2}}\right)}\\
&\geq&\int_{z-\E y_t}^{-z_0-\E y_t}\frac{-u\d u}{2|u|^{2\bar{\epsilon}\bar{\delta}}}\\
&=&\frac{1}{4(1-\bar{\epsilon}\bar{\delta})}\left(\left(-z+\E y_t\right)^{2(1-\bar{\epsilon}\bar{\delta})}-\left(z_0+\E y_t\right)^{2(1-\bar{\epsilon}\bar{\delta})}\right)\\
&=&\frac{1}{4(1-\bar{\epsilon}\bar{\delta})}\left(\left|z-\E y_t\right|^{2(1-\bar{\epsilon}\bar{\delta})}-\left|\mathrm{sgn}z\cdot z_0-\E y_t\right|^{2(1-\bar{\epsilon}\bar{\delta})}\right).
\eeqnn
Along the same lines as above, we can easily check the case $z\geq z_0$.
This completes our proof.
\fin

Notice that $z_t=\sigma_tu_x(t,\eta_t)$ and then $D^H_uz_t=\sigma_tu_{xx}(t,\eta_t)\sigma_\cdot$,
following exactly the same line as the proof of Theorem \ref{Theorem4.1} then yields a result for $z_t$, which we state as follows
\btheorem\label{Theorem4.2}
Let $t\in(0,T]$. Suppose that $0<\underline{u_x(t,\cdot)}<+\infty,\overline{u_{xx}(t,\cdot)}<+\infty$ and there exist positive constants $L,\lambda$ satisfying $u_{xx}(t,\cdot)\geq\frac{1}{L(1+|\cdot|^\lambda)}$.
Then the law of $z_t$ has a density $\rho_{z_t}$,
and for any $\epsilon,\delta>0$ there exists positive constants $C_1$ and $C_2$ such that
\beqnn
&&\frac{\E |z_t-\E z_t|}{C_1\sigma^2_t\varrho_t\left(1+|z|^{\bar{\epsilon}\bar{\delta}}\right)
\left(1+|z|^{\bar{\epsilon}\bar{\delta}}+\varrho_t^{\frac{\bar{\epsilon}}{2}}\right)}
\exp\left(-\frac{1}{C_2\sigma^2_t\varrho_t}\int_0^{z-\E z_t}u\left(1+|u+\E{z_t}|^{2\lambda\bar{\delta}}\right)\d u\right)\leq
\rho_{z_t}(z)\\
&\leq&\frac{\E |z_t-\E z_t|}{C_2\sigma^2_t\varrho_t}
\left(1+|z|^{2\lambda\bar{\delta}}\right)\exp\left(-\frac{1}{C_1\sigma^2_t\varrho_t}\int_0^{z-\E z_t}
\frac{u\d u}{\left(1+|u+\E{z_t}|^{\bar{\epsilon}\bar{\delta}}\right)\left(1+|u+\E{z_t}|^{\bar{\epsilon}\bar{\delta}}+\varrho_t^{\frac{\bar{\epsilon}}{2}}\right)}\right).
\eeqnn
\etheorem

Similar to Remark \ref{Remark4.1} and  Corollary \ref{Corollary4.1}, we have the following estimates.
\bcorollary\label{Corollary4.2}
With the same preamble as in Theorem \ref{Theorem4.2}, then the density $\rho_{z_t}$
fulfills the following bounds
\beqnn
\rho_{z_t}(z)\geq
\frac{\E |z_t-\E z_t|}{C_1\sigma^2_t\varrho_t\left(1+|z|^{\bar{\epsilon}\bar{\delta}}\right)
\left(1+|z|^{\bar{\epsilon}\bar{\delta}}+\varrho_t^{\frac{\bar{\epsilon}}{2}}\right)}
\exp\left(-\frac{\chi(z,\E z_t)}{C_2\sigma^2_t\varrho_t}\right),\ \ z\in\R
\eeqnn
and
\beqnn
\rho_{z_t}(z)
\leq\frac{\E |z_t-\E z_t|}{C_2\sigma^2_t\varrho_t}
\left(1+|z|^{2\lambda\bar{\delta}}\right)\exp\left(-\frac{\left|z-\E z_t\right|^{2(1-\bar{\epsilon}\bar{\delta})}-\left|\mathrm{sgn}z\cdot z_0-\E z_t\right|^{2(1-\bar{\epsilon}\bar{\delta})}}{4(1-\bar{\epsilon}\bar{\delta})C_1\sigma^2_t\varrho_t}\right), \ \ |z|>z_0
\eeqnn
with some positive constant $z_0$.
\ecorollary

\bremark\label{Remark 0}
(i) Comparing to the relevant results on BSDE driven by the standard Brownian motion ($H=\frac{1}{2}$)
proved in \cite[Theorem 5.6]{Mastrolia&Possama&Reveillac16a}, it is clear to see that 
our results apply to more general BSDEs since we replace $B^{\frac{1}{2}}_t$ with $\eta_t=\eta_0+\int_0^tb_s\d s+\int_0^t\sigma_s\d B^H_s$
and treat the case of fractional Brownian motion with arbitrary $H\in(\frac{1}{2},1)$ as driving process. Furthermore, an explicit low bound for the density of the solution $y_t$ without any restriction for the variable $z$ is shown in our Remark \ref{Remark4.1}.

(ii) The advantage of the above method of estimating densities allows us to obtain non-Gaussian type lower and upper bounds. The drawback is that Theorem \ref{Theorem4.1} and Theorem \ref{Theorem4.2} must be completed by an analysis of the following PDE
\begin{equation}\nonumber
\left\{
\begin{array}{ll}
u_t(t,x)+\frac{1}{2}\varrho_t'u_{xx}(t,x)+b(t)u_x(t,x)+f(t,x,u(t,x),\sigma_tu_x(t,x))=0
,\\
u(T,x)=h(x),
\end{array} \right.
\end{equation}
which is studied in \cite{Hu&cone&Song12a}. In the next two sections, we will present the Gaussian type densities estimates results where the only assumptions are those put on the 
data of BSDEs. 

\eremark

\section{Linear BSDEs driven by fractional Brownian motions}

\setcounter{equation}{0}

In the previous section, we have presented the non-Gaussian type densities estimates for the fractional BSDE \eqref{3.1},
which are indeed obtained by using the relation between this type of equation and a quasilinear PDE of mixed type.
While the BSDE \eqref{3.1} is linear, we are able to prove the Gaussian type bounds and the tail probabilities
with a direct and simpler method based on the explicit expression of the solutions in terms of the quasi-conditional expectation.

Consider the following linear BSDE
\begin{equation}\label{3.1'}
\left\{
\begin{array}{ll}
\d y_t=-\left(\alpha_t+\beta_ty_t+\gamma_tz_t\right)\d t-z_t\d B^H_t,\\
y_T=\xi,
\end{array} \right.
\end{equation}
where $\alpha_t,\beta_t,\gamma_t$ are given as continuous and adapted processes.

Notice that BSDE \eqref{3.1'} admits a unique solution under the condition \eqref{add000} below.
More specifically,
set $\varsigma_t:=\gamma_t+\int_0^tD_t^H\beta_s\d s,\bar{B}_t^H:=B_t^H+\int_0^t\varsigma_s\d s$
and
\beqnn
R_t:=\exp\left[-\int_0^t\left(K_H^{-1}\int_0^\cdot\varsigma_r\d r\right)(s)\d W_s-\frac{1}{2}\int_0^t\left(K_H^{-1}\int_0^\cdot\varsigma_r\d r \right)^2(s)\d s\right].
\eeqnn
It follows from the Novikov condition, i.e.
\beqlb\label{add000}
\E\exp \left[\frac{1}{2}\int_0^T\left(K_H^{-1}\int_0^\cdot\varsigma_r\d r \right)^2(s)\d s\right]<\infty,
\eeqlb
that $(R_t)_{t\in[0,T]}$ is an exponential martingale and then the Girsanov theorem implies that
$(\bar{B}^H)_{t\in[0,T]}$ is a fractional Brownian motion under the probability measure $\mathbb{Q}:=R_T\d \P$.
Let $\rho_t:=\exp\left\{\int_0^t\beta_s\d s\right\}$.
Applying the fractional integration by parts formula to $\rho_ty_t$ yields that
\beqlb\label{add3.1}
\d (\rho_t y_t)=-\alpha_t\rho_t\d t-\rho_t z_t\d \bar{B}^H_t.
\eeqlb
Then, there is a unique solution for BSDE \eqref{3.1'}, and moreover
\beqlb\label{3.2}
y_t=\rho_t^{-1}\hat{\mathbb{E}}^\mathbb{Q}\left[\rho_T\xi+\int_t^T\alpha_s\rho_s\d s|\mathcal{F}_t\right],
\eeqlb
where $\hat{\mathbb{E}}$ stands for the quasi-conditional expectation.
More details can be found in \cite[Theorem 5.1]{Hu&Peng09a} or \cite{Zhang14}.

Next, we want to show the existence of densities for the marginal laws of the solution $(y,z)$ to BSDE \eqref{add000} and then to derive the Gaussian bounds for them via the relation \eqref{3.2} and Proposition \ref{Theor2.2}. To this end, let $\alpha_t,\beta_t,\gamma_t$ be given continuous and deterministic functions and $\xi=h(\eta_T)$, in which $\eta$ is defined in the previous  section. Put
\beqnn
p_t(x):=\frac{1}{\sqrt{2\pi t}}\e^{-\frac{x^2}{2t}},
\eeqnn
and further denote
\beqnn
P_t g(x):=\int_\R p_t(x-y)g(y)\d y.
\eeqnn

We first state the following useful lemma concerning the representation of the quasi-conditional expectation.
\blemma\label{Lemma3.1}
Assume that $g:\mathbb{R}\to\mathbb{R}$ is a measurable function of polynomial growth, then the  following holds
\beqnn
\hat{\mathbb{E}}\left(g(\eta_T)|\mathcal{F}_t\right)=\P_{\|\sigma\|_T^2-\|\sigma\|_t^2}g\left(\eta_0+\int_0^Tb_s\d s+\int_0^t\sigma_s\d B^H_s\right).
\eeqnn
\elemma

\emph{Proof.}
Though the proof is similar to the one proposed in \cite[Theorem 3.8]{Hu&Peng09a},
yet we give a justification for the convenience of the reader.
For $t\in[0,T]$, let $\tilde{\eta}_t:=\eta_0+\int_0^Tb_s\d s+\int_0^t\sigma_s\d B^H_s$.
Applying the It\^{o} formula (see \cite[Theorem 2.3]{Hu&Peng09a} or \cite[Corollary 35]{Maticiuc&Nie15a})
to $P_{\|\sigma\|_t^2-\|\sigma\|_s^2}g(\tilde{\eta}_s)$, we get
\beqlb\label{Lemma3.1-1}
g(\tilde{\eta}_t)=P_{\|\sigma\|_t^2}g(\tilde{\eta}_0)+\int_0^t\frac{\partial}{\partial x}P_{\|\sigma\|_t^2-\|\sigma\|_s^2}g(\tilde{\eta}_s)\sigma_s\d B^H_s.
\eeqlb
Choosing $t=T$ in the above equation and then taking the quasi-conditional expectation with respect to $\mathcal{F}_t$, we have
\beqlb\label{Lemma3.1-2}
\hat{\mathbb{E}}\left(g(\tilde{\eta}_T)|\mathcal{F}_t\right)=P_{\|\sigma\|_T^2}g(\tilde{\eta}_0)
+\int_0^t\frac{\partial}{\partial x}P_{\|\sigma\|_T^2-\|\sigma\|_s^2}g(\tilde{\eta}_s)\sigma_s\d B^H_s.
\eeqlb
By the semigroup property of $P_t$, it is easy to verify that, for $0\leq s\leq t\leq T$,
\beqnn
\frac{\partial}{\partial x}P_{\|\sigma\|_T^2-\|\sigma\|_s^2}g(x)=
P_{\|\sigma\|_T^2-\|\sigma\|_t^2}\frac{\partial}{\partial x}P_{\|\sigma\|_t^2-\|\sigma\|_s^2}g(x).
\eeqnn
Hence, this, together with \eqref{Lemma3.1-2} and \eqref{Lemma3.1-1}, yields the desired result.
\fin

Frow now on, let us suppose the following \\
(H2) $h:\mathbb{R}\to\mathbb{R}$ is twice differentiable and $0<c\leq h'\leq C, 0<\tilde{c}\leq h''\leq\tilde{C}$, where $c, C, \tilde{c}$ and $\tilde{C}$ are constants.\\
Besides, we set
\beqnn
\vartheta_1(t):=\varrho_t\exp\left[2(T-t)\inf\limits_{s\in[0,T]}\beta_s\right]
\eeqnn
and
\beqnn
\vartheta_2(t):=\varrho_t\exp\left[2(T-t)\sup\limits_{s\in[0,T]}\beta_s\right].
\eeqnn
Recall that $\varrho_t$ is defined in Section 3.

We are now in the position to state our main result of this section.
\btheorem\label{Theorem3.1}
Assume that (H2) holds.
Then, for each $t\in(0,T],y_t$ and $z_t$ possess densities $p_{y_t}$ and $p_{z_t}$, respectively.
Moreover, for almost all $x\in\R, p_{y_t}$ and $p_{z_t}$ satisfy, respectively, the following bounds
\beqnn
\frac{\E|y_t-\E y_t|}{2C^2\vartheta_2(t)}\exp\left(-\frac{(x-\E y_t)^2}{2c^2\vartheta_1(t)}\right)
\leq p_{y_t}(x)\leq
\frac{\E|y_t-\E y_t|}{2c^2\vartheta_1(t)}\exp\left(-\frac{(x-\E y_t)^2}{2C^2\vartheta_2(t)}\right).
\eeqnn
and
\beqnn
\frac{\E|z_t-\E z_t|}{2\tilde{C}^2\vartheta_2(t)\sigma^2_t}\exp\left(-\frac{(x-\E z_t)^2}{2\tilde{c}^2\vartheta_1(t)\sigma^2_t}\right)
\leq p_{z_t}(x)\leq
\frac{\E|z_t-\E z_t|}{2\tilde{c}^2\vartheta_1(t)\sigma^2_t}\exp\left(-\frac{(x-\E z_t)^2}{2\tilde{C}^2\vartheta_2(t)\sigma^2_t}\right).
\eeqnn
\etheorem

\emph{Proof.}
By \eqref{3.2}, we obtain
$y_t=\e^{\int_t^T\beta_s\d s}\hat{\mathbb{E}}^\mathbb{Q}[\xi|\mathcal{F}_t]+\rho_t^{-1}\int_t^T\alpha_s\rho_s\d s$.
Hence, we have
\beqlb\label{Theorem3.1-1}
D_u^Hy_t=\e^{\int_t^T\beta_s\d s}D_u^H\left(\hat{\mathbb{E}}^\mathbb{Q}[\xi|\mathcal{F}_t]\right).
\eeqlb
On the other hand, from Lemma \ref{Lemma3.1}, we get
\beqlb\label{add0''-Theorem3.1}
\hat{\mathbb{E}}^\mathbb{Q}[\xi|\mathcal{F}_t]
=\hat{\mathbb{E}}^\mathbb{Q}[h(\eta_T)|\mathcal{F}_t]
&=&\hat{\mathbb{E}}^\mathbb{Q}\left[h\left(\eta_0+\int_0^Tb_s\d s-\int_0^T\sigma_s\gamma_s\d s+\int_0^T\sigma_s\d\bar{B}^H_s\right)|\mathcal{F}_t\right]\nonumber\\
&=&P_{\|\sigma\|_T^2-\|\sigma\|_t^2}h\left(\eta_0+\int_0^Tb_s\d s-\int_0^T\sigma_s\gamma_s\d s+\int_0^t\sigma_s\d \bar{B}^H_s\right)\nonumber\\
&=&P_{\|\sigma\|_T^2-\|\sigma\|_t^2}h\left(\eta_0+\int_0^Tb_s\d s-\int_t^T\sigma_s\gamma_s\d s+\int_0^t\sigma_s\d B^H_s\right).
\eeqlb
Then, for $u\in[0,t]$
\beqlb\label{add0'-Theorem3.1}
D_u^H\left(\hat{\mathbb{E}}^\mathbb{Q}[\xi|\mathcal{F}_t]\right)
=\sigma_u P_{\|\sigma\|_T^2-\|\sigma\|_t^2}h'\left(\eta_0+\int_0^Tb_s\d s-\int_t^T\sigma_s\gamma_s\d s+\int_0^t\sigma_s\d B^H_s\right).
\eeqlb
This allows us to deduce from \eqref{Theorem3.1-1} that, for $u\in[0,t]$
\beqlb\label{add1-Theorem3.1}
D_u^Hy_t=\sigma_u\e^{\int_t^T\beta_s\d s}P_{\|\sigma\|_T^2-\|\sigma\|_t^2}h'\left(\eta_0+\int_0^Tb_s\d s-\int_t^T\sigma_s\gamma_s\d s+\int_0^t\sigma_s\d B^H_s\right).
\eeqlb
Notice that
\beqnn
\Phi_{y_t}(B^H)=D^Hy_t=\sigma_\cdot\e^{\int_t^T\beta_s\d s}P_{\|\sigma\|_T^2-\|\sigma\|_t^2}h'\left(\eta_0+\int_0^Tb_s\d s-\int_t^T\sigma_s\gamma_s\d s+\int_0^t\sigma_s\d B^H_s\right),
\eeqnn
then
\beqnn
&&\widetilde{\Phi_{y_t}^\theta}(B^H)\\
&=&\Phi_{y_t}(\e^{-\theta}B^H+\sqrt{1-\e^{-2\theta}}B'^H)\\
&=&\sigma_\cdot\e^{\int_t^T\beta_s\d s}P_{\|\sigma\|_T^2-\|\sigma\|_t^2}h'\left(\eta_0+\int_0^Tb_s\d s-\int_t^T\sigma_s\gamma_s\d s
+\e^{-\theta}\int_0^t\sigma_s\d B^H_s+\sqrt{1-\e^{-2\theta}}\int_0^t\sigma_s\d B'^H_s\right).
\eeqnn
Thus, according to \eqref{2.1.1}, we have
\beqnn
\langle\Phi_{y_t}(B^H),\widetilde{\Phi_{y_t}^\theta}(B^H)\rangle_\mathcal{H}=\varrho_t\kappa(t,\theta),
\eeqnn
where
\beqnn
&&\kappa(t,\theta)\\
&=&\e^{2\int_t^T\beta_s\d s}
P_{\|\sigma\|_T^2-\|\sigma\|_t^2}h'\left(\eta_0+\int_0^Tb_s\d s-\int_t^T\sigma_s\gamma_s\d s+\int_0^t\sigma_s\d B^H_s\right)\\
&&\times P_{\|\sigma\|_T^2-\|\sigma\|_t^2}h'\left(\eta_0+\int_0^Tb_s\d s-\int_t^T\sigma_s\gamma_s\d s
+\e^{-\theta}\int_0^t\sigma_s\d B^H_s+\sqrt{1-\e^{-2\theta}}\int_0^t\sigma_s\d B'^H_s\right)
\eeqnn
with $c^2\e^{2(T-t)\inf_{s\in[0,T]}\beta_s}\leq\kappa(t,\theta)\leq C^2\e^{2(T-t)\sup_{s\in[0,T]}\beta_s}$ due to (H2).
Consequently, we arrive at the following bound
\beqlb\label{add2-Theorem3.1}
c^2\vartheta_1(t)\leq g_{y_t}\leq C^2\vartheta_2(t).
\eeqlb

Next, we devote to computing $D^H_\cdot z_t$ and then estimating $g_{z_t}$.
From \eqref{add3.1}, we know
\beqlb\label{Theorem3.1-2}
y_0-\rho_T\xi-\int_0^T\alpha_t\rho_t\d t=\int_0^T\rho_t z_t\d\bar{B}_t^H.
\eeqlb
On the other hand, by the fractional Clark-Ocone formula (see \cite{Hu05a} and \cite{Hu&Oksendal03a}) on (${\bar{B}^H,\mathbb{Q}}$), we can write
\beqlb\label{Theorem3.1-3}
y_0-\rho_T\xi-\int_0^T\alpha_t\rho_t\d t
&=&\int_0^T\hat{\mathbb{E}}^\mathbb{Q}\left[D^H_t\left(y_0-\rho_T\xi-\int_0^T\alpha_t\rho_t\d t\right)|\mathcal{F}_t\right]\d\bar{B}_t^H\nonumber\\
&=&-\int_0^T\rho_T\hat{\mathbb{E}}^\mathbb{Q}\left[D^H_t\xi|\mathcal{F}_t\right]\d\bar{B}_t^H.
\eeqlb
Then, combining \eqref{Theorem3.1-2} with \eqref{Theorem3.1-3} yields the following
\beqnn
z_t=-\frac{\rho_T}{\rho_t}\hat{\mathbb{E}}^\mathbb{Q}\left[D^H_t\xi|\mathcal{F}_t\right]
=-\sigma_t\e^{\int_t^T\beta_s\d s}\hat{\mathbb{E}}^\mathbb{Q}\left[h'(\eta_T)|\mathcal{F}_t\right].
\eeqnn
Similar to \eqref{add1-Theorem3.1} and \eqref{add2-Theorem3.1}, we obtain
\beqlb\label{Theorem3.1-4}
D_u^Hz_t=-\sigma_u\sigma_t\e^{\int_t^T\beta_s\d s}P_{\|\sigma\|_T^2-\|\sigma\|_t^2}h''\left(\eta_0+\int_0^Tb_s\d s-\int_t^T\sigma_s\gamma_s\d s+\int_0^t\sigma_s\d B^H_s\right),
\eeqlb
and
\beqlb\label{Theorem3.1-5}
\tilde{c}^2\vartheta_1(t)\sigma^2_t\leq g_{z_t}\leq\tilde{C}^2\vartheta_2(t)\sigma^2_t.
\eeqlb
Finally, applying Proposition \ref{Theor2.2} to \eqref{add2-Theorem3.1} and \eqref{Theorem3.1-5} respectively, we end up with the desired results and the proof is complete.
\fin

\bremark\label{Remark3.1}
We can alternatively derive \eqref{add0'-Theorem3.1} in the above proof, by the following
\beqnn
D_u^H\left(\hat{\mathbb{E}}^\mathbb{Q}[h(\eta_T)|\mathcal{F}_t]\right)
&=&\hat{\mathbb{E}}^\mathbb{Q}[D_u^Hh(\eta_T)|\mathcal{F}_t]
=\sigma_u\mathrm{I}_{[0,t]}(u)\hat{\mathbb{E}}^\mathbb{Q}[h'(\eta_T)|\mathcal{F}_t]\\
&=&\sigma_u\mathrm{I}_{[0,t]}(u)P_{\|\sigma\|_T^2-\|\sigma\|_t^2}h'\left(\eta_0+\int_0^Tb_s\d s-\int_t^T\sigma_s\gamma_s\d s+\int_0^t\sigma_s\d B^H_s\right).
\eeqnn
where the last equality is similar to \eqref{add0''-Theorem3.1}.
\eremark

In view of the proof of Theorem \ref{Theorem3.1}, one can derive the following tail estimates for the probability laws of $y_t$ and $z_t$.

\bcorollary\label{Corollary 3.1}
Suppose (H2).
Then there hold, for all $x>0$,
\beqlb\label{Coro3.1-1}
\P(y_t-\E y_t\geq x)\leq\exp\left(-\frac{x^2}{2C^2\vartheta_2(t)}\right), \ \ \P(y_t-\E y_t\leq-x)\leq\exp\left(-\frac{x^2}{2C^2\vartheta_2(t)}\right),
\eeqlb
and
\beqlb\label{Coro3.1-2}
\P(z_t-\E z_t\geq x)\leq\exp\left(-\frac{x^2}{2\tilde{C}^2\vartheta_2(t)\sigma^2_t}\right),\ \
\P(z_t-\E z_t\leq-x)\leq\exp\left(-\frac{x^2}{2\tilde{C}^2\vartheta_2(t)\sigma^2_t}\right).
\eeqlb
\ecorollary

\emph{Proof.}
Noticing first that $g_F(x)=g_{F-\E F}(x)$.
Then the relation \eqref{Coro3.1-1} follows by \eqref{add2-Theorem3.1} and Proposition \ref{Proposition 2.1} with $a_1=0$ and $a_2=C^2\vartheta_2(t)$.
\eqref{Coro3.1-2} can be verified  similarly.
\fin

\section{BSDEs driven by Gaussian processes}

\setcounter{equation}{0}

In this section, we consider the following BSDE driven by a centered Gaussian process $\{X_t\}_{t\in[0,T]}$
\begin{equation}\label{5.1}
\left\{
\begin{array}{ll}
\d y_t=-f(t,X_t,y_t,z_t)\d V(t)+z_t\d^\diamond X_t,\\
y_T=h(X_T),
\end{array} \right.
\end{equation}
where $\mathrm{Var}X_t=V(t), t\in[0,T]$, is a strictly increasing, continuous function with $V(0)=0$, the stochastic integral is the Wick-It\^{o} integral defined by the $S$-transformation and the Wick product.

As stated in \cite[Section 4.4]{Bender}, BSDE \eqref{5.1} covers a wide class of Gaussian processes as driving processes including fractional Brownian motion with $H\in(0,1)$ and fractional Wiener integral and so on, wherein a comparison with BSDE \eqref{3.1} is also
given. Moreover, we notice that Bender \cite{Bender} shows the existence and uniqueness results for BSDE \eqref{5.1} under the Lipschitz or even superquadratic growth conditions via the transfer theorem which can transfer the concerned problems to an auxiliary BSDE driven by a Brownian motion. Recall that the auxiliary BSDE is of the following form
\begin{equation}\label{5.2}
\left\{
\begin{array}{ll}
\d \bar{y}_t=-f(U(t),\bar{W}_t,\bar{y}_t,\bar{z}_t)\d t+\bar{z}_t\d\bar{W}_t,\\
\bar{y}_{V(T)}=h(\bar{W}_{V(T)}),
\end{array} \right.
\end{equation}
where $\{\bar{W}_t\}_{t\in[0,V(T)]}$ is a standard Brownian motion,
$U(t), t\in[0,V(T)]$, is the inverse of $V$ defined as
\beqnn
U(t):=\inf\{s\geq0: V(s)\geq t\}, \, t\in[0,V(T)].
\eeqnn

In this part, we aim to investigate the existence of densities and then derive their Gaussian  estimates for the marginal laws of the solution $(y,z)$ to BSDE \eqref{5.1}.
For this, we start by adopting the following set of conditions from \cite{Aboura&Bourguin13a,Bender}: (H3)
\begin{itemize}
\item[(i)] $h\in C_b^2(\R)$ and $\inf\limits_{x\in\R} h''(x)>0$;

\par

\item[(ii)]
$\E\int_0^T|f(t,X_t,0,0)|^2\d V(t)<\infty$, and
for each $t\in[0,T], f(t,\cdot,\cdot,\cdot)\in C_b^2(\R^3)$ with all
positive derivatives.
\end{itemize}

We state our final main result as follows

\btheorem\label{Theorem 5.1}
Assume (H3) holds true.
Then, for all $t\in[0,T]$,
\begin{itemize}
\item[(1)] $y_t$ has a density $\rho_{y_t}$, and there exist two strictly positive constants $c_1<c_2$ such that for any $z\in\R$,
\beqlb\label{Theorem 5.1-1'}
\frac{\E|y_t-\E y_t|}{c_2V(t)}\exp\left(-\frac{(z-\E y_t)^2}{c_1V(t)}\right)
&\leq&\rho_{y_t}(z)\nonumber\\
&\leq&\frac{\E|y_t-\E y_t|}{c_1V(t)}\exp\left(-\frac{(z-\E y_t)^2}{c_2V(t)}\right),
\eeqlb

\par

\item[(2)] if the generator $f(t,x,y,z)$ has a linear dependence on the $z$ component,
$z_t$ possesses a density $\rho_{z_t}$, and furthermore if $f$ only depends on the components $(t,y)$, there exist two strictly positive constants $c_3<c_4$ such that for any $z\in\R$,
\beqlb\label{Theorem 5.1-1''}
\frac{\E|z_t-\E z_t|}{c_4V(t)}\exp\left(-\frac{(z-\E z_t)^2}{c_3V(t)}\right)
&\leq&\rho_{z_t}(z)\nonumber\\
&\leq&\frac{\E|z_t-\E z_t|}{c_3V(t)}\exp\left(-\frac{(z-\E z_t)^2}{c_4V(t)}\right).
\eeqlb
\end{itemize}
\etheorem

\emph{Proof.}
By (H3), it follows by \cite[Theorem 4.1]{Karoui&Peng&Quenez97a} that the auxiliary BSDE \eqref{5.2} admits a unique solution $(\bar{y},\bar{z})$
which has the following representation
\beqlb\label{Theorem 5.1-1}
\bar{y}_t=\phi(t,\bar{W}_t),\ \ \bar{z}_t=\psi(t,\bar{W}_t),\ \ t\in[0,V(T)],
\eeqlb
where both $\phi,\psi: [0,V(T)]\times\R\rightarrow\R$ are deterministic functions.\\
Then, in light of \cite[Theorem 4.2 and Theorem 4.4]{Bender} or the transfer theorem asserted in \cite[Theorem 3.1]{Bender}, we conclude that
BSDE \eqref{5.2} has a unique solution $(y,z)$ which can be written as
\beqlb\label{Theorem 5.1-2}
y_t=\phi(V(t),X_t),\ \ z_t=\psi(V(t),X_t),\ \ t\in[0,T].
\eeqlb
Therefore, taking into account the fact that $\mathrm{Var}X_t=V(t)$ we deduce that the law of $y_t$ (resp. $z_t$) is the same as that of $\bar{y}_{V(t)}$ (resp. $\bar{z}_{V(t)}$).

On the other hand, by \cite[Theorem 3.3]{Aboura&Bourguin13a} it follows that $\bar{y}_{V(t)}$
possesses a density $\rho_{\bar{y}_{V(t)}}$, and moreover there exist some strictly positive constants $c_1<c_2$
such that for all $z\in\R$,
\beqlb\label{Theorem 5.1-3}
&& \frac{\E|\bar{y}_{V(t)}-\E\bar{y}_{V(t)}|}{c_2V(t)}\exp\left(-\frac{(z-\E\bar{y}_{V(t)})^2}{c_1V(t)}
\right)
\leq \rho_{\bar{y}_{V(t)}}(z)\nonumber\\
&\leq&\frac{\E|\bar{y}_{V(t)}-\E\bar{y}_{V(t)}|}{c_1V(t)}\exp\left(-\frac{(z-\E\bar{y}_{V(t)})^2}{c_2V(t)}\right),
\eeqlb
which then yields our first claim.\\
As for $z_t$, provided that the generator $f$ has a linear dependence on the $z$ component,
owing to \cite[Theorem 4.3]{Aboura&Bourguin13a}, we conclude that $\bar{z}_{V(t)}$ has a law which is absolutely continuous with respect to the Lebesgue measure.
Moreover, if $f$ depends only on $(t,y)$ components, then applying \cite[Theorem 4.6]{Aboura&Bourguin13a}, we obtain the following Gaussian bounds for the density of $\bar{z}_{V(t)}$
\beqlb\label{Theorem 5.1-4}
\frac{\E|\bar{z}_{V(t)}-\E\bar{z}_{V(t)}|}{c_4V(t)}\exp\left(-\frac{(z-\E\bar{z}_{V(t)})^2}{c_3V(t)}\right)
&\leq&\rho_{\bar{z}_{V(t)}}(z)\nonumber\\
&\leq&\frac{\E|\bar{z}_{V(t)}-\E\bar{z}_{V(t)}|}{c_3V(t)}\exp\left(-\frac{(z-\E\bar{z}_{V(t)})^2}{c_4V(t)}\right),
\eeqlb
where $c_3<c_4$ are two strictly positive constants.\\
We therefore obtain the other assertion.
\fin

\bremark\label{Remark 5.2}
(i) If we choose Brownian motion as the driving Gaussian process, namely $X_t=B^{\frac{1}{2}}$,
then our estimates \eqref{Theorem 5.1-1'} and \eqref{Theorem 5.1-1''} coincide with the estimates of
\cite[Theorem 3.3 and Theorem 4.6]{Aboura&Bourguin13a}.
Note that the equation \eqref{5.1} that we considered is more general than that of \cite{Aboura&Bourguin13a}
since we allow $X$ to be a wide class of Gaussian processes which includes fractional Brownian motion with arbitrary $H\in(0,1)$ as special case.
Therefore, the results stated in Theorem \ref{Theorem 5.1} cover that of \cite{Aboura&Bourguin13a}.

(ii) If we take $X_t=\eta_t$,
then the equation \eqref{5.1} is of the same form as \eqref{3.1} and \eqref{3.1'},
which we considered in Section 3 and 4, respectively.
By simple calculus we know that the result of Theorem \ref{Theorem4.1} combined with Remark \ref{Remark4.1} and Corollary \ref{Corollary4.1}
is more elaborate than that of Theorem \ref{Theorem 5.1}.
Similarly, we also note that the estimates in Theorem \ref{Theorem3.1} are better than the estimates \eqref{Theorem 5.1-1'} and \eqref{Theorem 5.1-1''} in Theorem \ref{Theorem 5.1}.
Indeed, the derivation of Theorem \ref{Theorem 5.1} is mainly based on the transfer theorem
which is a time change type transformation allowing us to represent the equation \eqref{5.1} in terms of a class of BSDE driven by Brownian motion,                while the arguments used in Section 3 and 4 focus on the structures of the original equations.

\eremark

\textbf{Acknowledgements}
The research of X. Fan was supported in part by the National Natural Science Foundation of China (Grant No. 11501009, 11371029), the Natural Science Foundation of Anhui Province (Grant No.
1508085QA03), the Distinguished Young Scholars Foundation of Anhui Province (Grant No. 1608085J06).

\end{document}